\documentclass[12pt]{amsart}
\usepackage{a4wide}
\usepackage[dvips]{graphicx}
\usepackage{wasysym}
\usepackage{amsfonts}
\usepackage{amssymb}
\usepackage{epsfig}
\usepackage[all]{xy}
\newtheorem{defi}{\textbf{Definition}}[section]
\newtheorem{definition}[defi]{\textbf{Definition}}
\newtheorem{theorem}[defi]{\textbf{Theorem}}
\newtheorem{lemma}[defi]{\textbf{Lemma}}
\newtheorem{proposition}[defi]{\textbf{Proposition}}
\newtheorem{corollary}[defi]{\textbf{Corollary}}
\newtheorem{assumption}[defi]{\textbf{Assumption}}

\newtheorem{remark}[defi]{\textbf{Remark}}

\def\Z {\mathbb Z}
\def\N {\mathbb N}
\def\mc {\mathcal}

\begin{document}

\title[Free splittings in relatively hyperbolic groups]{Detecting free splittings in relatively hyperbolic groups}

\author[Fran\c{c}ois Dahmani]{Fran\c{c}ois Dahmani}
\address{Fran\c{c}ois Dahmani\\
Laboratoire E. Picard\\
Universit\'e Paul Sabatier\\
F-31062 Toulouse, France}
\email{dahmani@picard.ups-tlse.fr}

\author[Daniel Groves]{Daniel Groves}
\address{Daniel Groves\\
Department of Mathematics \\
California Institute of Technology \\
Pasadena, CA, 91125, USA} 
\email{groves@caltech.edu} 

\date{1 December, 2006}

\thanks{The first author acknowledges support from the ANR grant 06-JCJC-0099-01. The second author's work was supported in part by NSF Grant DMS-0504251.}

\begin{abstract}
We describe an algorithm which determines whether or not a group which
is hyperbolic relative to abelian groups admits a nontrivial splitting over a finite group.
\end{abstract}

\maketitle

\section{Introduction}

In this paper we are concerned with the problem of determining whether or not
a group given by some finite presentation admits a nontrivial free product 
decomposition (or a nontrivial splitting over a finite subgroup).  Of course this problem has no solutions for arbitrary
finite group presentations.  For example, the group $(\Z / 2\Z) \ast G$ is freely 
indecomposable if and only if $G \cong 1$, but it is impossible to decide whether or
not a finitely presented group $G$ is trivial.  However, in the presence of some geometry, a
positive solution is sometimes possible.  
        One of the most important results in this direction is due to V. Gerasimov \cite{Gera}. In the late 90's, he proved:

\begin{theorem}\label{theo;Gera} (V. Gerasimov)
 There is an algorithm that, given a finite presentation of a hyperbolic group, computes the number of ends of the group.
\end{theorem}

In this paper, we consider {\em relatively hyperbolic groups}. (See Section 
\ref{s:Prelim} for a definition.)  Following Gerasimov's ideas, our main result is:

\begin{theorem} \label{t:ExistsSplit}
There exists an algorithm which takes as input a finite presentation of a group
$\Gamma$ which is hyperbolic relative to abelian subgroups, and outputs `Yes' or `No',
according to whether or not $\Gamma$ admits a nontrivial splitting over a finite group.
\end{theorem}

According to Dunwoody's accessibility, for a given finitely presented groups $G$ 
there is a bound on the size of a graph of groups with finite edge groups and
fundamental group $G$.  Thus, there is a {\em Dunwoody} decomposition of such
a $G$, which is the most refined splitting of $G$ as a graph of groups with finite edge groups. 
 Using Theorem \ref{t:ExistsSplit}, it is straightforward to deduce the following:

\begin{theorem} \label{t:Dunwoody}
There is an algorithm which takes as input a finite presentation of a group $\Gamma$
which is hyperbolic relative to abelian subgroups, and outputs a presentation for
a graph of groups encoding the Dunwoody decomposition of $\Gamma$, together
with an isomorphism to $\Gamma$.
\end{theorem}

If we consider only splittings over the trivial group, then there is the {\em Grushko
decomposition}, and it is possible to find this also.

\begin{theorem} \label{t:Grushko}
There is an algorithm which takes as input a finite presentation of a group $\Gamma$
which is hyperbolic relative to abelian subgroups, and outputs a presentation for
a graph of groups encoding the Grushko decomposition of $\Gamma$, together
with an isomorphism to $\Gamma$.
\end{theorem}

Theorem \ref{t:Grushko} is a major step in the authors' solution to the Isomorphism
Problem for toral relatively hyperbolic groups in \cite{DG}.

Let us briefly outline the method, and the background. Gerasimov's algorithm for hyperbolic groups is unfortunately still unpublished to this date, and, up to our knowledge, not publicly accessible as an e-print. The results in this paper imply Theorem \ref{theo;Gera}, and also extend it to a wide class of groups with the geometry of non-positive curvature. As in Gerasimov's approach, we use connectivity properties of boundaries.
 In \cite{Bowditch_rel_hyp}, Bowditch associates a natural boundary to every relatively hyperbolic group. This boundary is disconnected if the group admits a splitting over a finite subgroup, that is compatible with the relative structure ({\it i.e. } in which the parabolic subgroups are elliptic). When the parabolic subgroups are themselves one-ended, this latter requirement is always fulfilled, thus connectivity of the boundary is equivalent to one-endedness in this case.  In \cite{BestvinaMess}, Bestvina and Mess introduced a condition for (local) connectedness of the boundary of a hyperbolic group, provided this boundary has no global cut point (which Bowditch and Swarup latter proved to never happen for hyperbolic groups).  This condition can be stated in the relatively hyperbolic setting as well, but this time, global cut points can occur in the boundary. However, this situation was studied further by Bowditch, who proved,  in a quite wide generality, that  global cut points are so-called parabolic fixed points, and come from peripheral splittings (see Theorem \ref{theo;bowditch} below,  and references therein). Bowditch also proved an accessibility theorem for peripheral splittings of relatively hyperbolic groups. 

The main idea for detecting the number of ends of a relatively hyperbolic group is to decide whether its boundary is connected. For that we wish to  use the condition of  Bestvina and Mess, that can be expressed in terms  of a finite subset of the group (see Corollary \ref{coro;carac_ddagger}).  There is a one-sided algorithm which will terminate
if this condition is satisfied, and implies the boundary is connected (see Proposition
\ref{prop;check_ddagger}).  But if there are cut points in the boundary,  we want to compute a peripheral splitting so that no vertex group has this problem (by Bowditch's accessibility there is always such a splitting), and test the property vertex by vertex. A first algorithm  will enumerate  the (peripheral) splittings of the given group, while on each such splitting, another will check whether or not the splitting is nontrivial over a finite group, and if not, a third algorithm will check whether the condition of Bestvina and Mess is fulfilled for each vertex group. This latter algorithm might not terminate on every splitting, but it will terminate at least on a maximal peripheral splitting (the existence of which is guaranteed by Bowditch's accessibility result).

 In Section \ref{s:Prelim} we recall the definition
of relatively hyperbolic groups, and some of the geometric and algorithmic properties
which we require for this paper.  In Section \ref{s:Splittings} we recall some properties
about splittings of relatively hyperbolic groups, due to Bowditch 
 (mostly gathered in Theorem \ref{theo;bowditch}).  The key idea here
is that splittings of a relatively hyperbolic group can be detected in terms of cut points on the boundary.  
 In Section  \ref{s:BM} we introduce a relatively hyperbolic version of a 
condition of Bestvina and Mess \cite{BestvinaMess} which can be used, if the boundary has no cut point,  
 to discover whether it is connected (Proposition \ref{prop;check_ddagger}).  
Finally, in Section \ref{s:Algorithms} we present the algorithms which prove Theorems
\ref{t:ExistsSplit}, \ref{t:Dunwoody} and \ref{t:Grushko}.

We wish to thank the referee for precise and efficient remarks.  These improved
the paper, particularly the results in Subsection \ref{ss:HoroProps}.

\section{Preliminaries} \label{s:Prelim}

Relatively hyperbolic groups were originally defined by Gromov in his seminal
paper \cite{Gromov}.  An alternative definition was given by Farb \cite{Farb},
who further developed the theory and proved many fundamental results.  By
now there are many equivalent definitions and characterisations of relatively
hyperbolic groups, and there has been a great deal of recent activity.

We will provide a definition of relative hyperbolicity due to the second author and
Manning \cite{GrovesManning} which is suitable for our needs.

\subsection{Combinatorial horoballs and the cusped space}

In this paragraph, we recall a construction from \cite{GrovesManning}.  For
the majority of this paper, we
will only need the $1$-skeleton of the cusped space defined in Definition
\ref{d:X1}.

\begin{defi} \label{d:horoball}
Let $\Gamma$ be any $1$-complex.
The \emph{combinatorial horoball based on $\Gamma$}, denoted
${\mc H}(\Gamma)$, is the $2$-complex formed as follows:
\begin{itemize}
\item ${\mc H}^{(0)}= \Gamma^{(0)}\times \left( \{0\}\cup \N \right)$
\item ${\mc H}^{(1)}$ contains the following three types of edges.  The
  first two types are called \emph{horizontal}, and the last type is
  called \emph{vertical}.
\begin{enumerate}
\item If $e$ is an edge of $\Gamma$ joining $v$ to $w$ then there is a
  corresponding edge $\bar{e}$ connecting $(v,0)$ to $(w,0)$.
\item If $k>0$ and $0<d_{\Gamma}(v,w)\leq 2^k$, then there is a single edge
  connecting $(v,k)$ to $(w,k)$.
\item If $k\geq 0$ and $v\in \Gamma^{(0)}$, there is an edge  joining
  $(v,k)$ to $(v,k+1)$. 
\end{enumerate}
\item ${\mc H}^{(2)}$ contains three kinds of $2$-cells:
\begin{enumerate}
\item If $\gamma\subset {\mc H}^{(1)}$ is a circuit composed of three
  horizontal edges, then there is a $2$-cell (a \emph{horizontal
  triangle}) attached along $\gamma$.
\item If $\gamma\subset {\mc H}^{(1)}$ is a circuit composed of two
  horizontal edges and two vertical edges, then there is a $2$-cell (a
  \emph{vertical square}) attached along $\gamma$.
\item If $\gamma\subset {\mc H}^{(1)}$ is a circuit composed of three
  horizontal edges and two vertical ones, then there is a $2$-cell (a
  \emph{vertical pentagon}) attached along $\gamma$, unless $\gamma$
  is the boundary of the union of a vertical square and a horizontal
  triangle.
\end{enumerate}
\end{itemize}
\end{defi}

\begin{remark}\label{r:subset}
As the full subgraph of ${\mc H}(\Gamma)$ containing the vertices
$\Gamma^{(0)}\times\{0\}$ is isomorphic to $\Gamma$, we may think of
$\Gamma$ as a subset of ${\mc H}(\Gamma)$.
\end{remark}

We also define, for every $k \ge 0$, the {\em $k$-thick part} of a combinatorial
horoball, denoted ${\mc H}(\Gamma)_k$ to be the full subcomplex on the set of
those vertices of the form $(\gamma,j)$ with $0 \le j \le k$.

See \cite[Section 3]{GrovesManning} for a detailed discussion of the geometry
of combinatorial horoballs.  We recall only the properties which we will need
(these statements are \cite[Proposition 3.7, Theorem 3.8, Lemma 3.10]{GrovesManning}):

\begin{proposition} \label{p:HoroProps}
 \ 

\begin{enumerate}
\item Let $\Gamma$ be a connected $1$-complex so that no edge joins a vertex
to itself.  Then ${\mc H}(\Gamma)$ is simply-connected and satisfies a linear combinatorial isoperimetric inequality with constant at most $3$;
\item Let $\Gamma$ be any $1$-complex.  Then ${\mc H}(\Gamma)^{(1)}$ is $\Upsilon$-hyperbolic, and $\Upsilon$ is independent of $\Gamma$;
\item\label{item:geo} Let $\Gamma$ be a $1$-complex, and $a, b \in {\mc H}(\Gamma)^{(0)}$.  Then
there is a geodesic in ${\mc H}(\Gamma)^{(1)}$ joining $a$ to $b$ which consists
of at most two vertical segments, one going down (deeper in the horoball) and the other going up, and at most one horizontal segment of length at
most $3$.
\end{enumerate}
\end{proposition}

We now proceed to define the {\em cusped space}.

\begin{definition}\label{d:X1}
Let $G$ be a finitely generated group, let
${\bf P}=\{P_1,\ldots,P_n \}$ be a (finite) 
family of finitely generated subgroups of $G$, and let ${\mc A}$ be a
generating set for $G$ so that $P_i\cap {\mc A}$ generates $P_i$ for each
$i\in \{1,\ldots,n\}$.  
For each $i\in \{1,\ldots,n\}$, let $T_i$ be a left transversal for
$P_i$. 

For each $i$, and each $t\in T_i$, let $\Gamma_{i,t}$ be the full
subgraph of the Cayley graph $\Gamma(G,{\mc A})$ which contains $tP_i$.
Each $\Gamma_{i,t}$ is isomorphic to the Cayley graph of $P_i$ with
respect to the generators $P_i\cap {\mc A}$.
Then we define
\[ X = \Gamma(G,{\mc A}) \cup (\cup \{{\mc H}(\Gamma_{i,t}) \mid 1\leq i\leq n,
t\in T_i\}),\]
where the graphs $\Gamma_{i,t}\subset \Gamma(G,{\mc A})$ and
$\Gamma_{i,t}\subset {\mc H}(\Gamma_{i,t})$ are identified as suggested
in Remark \ref{r:subset}.
\end{definition}

We call $X$ the cusped space associated to $(G, {\bf P}, {\mc A})$. We also define the $k$-thick part of the cusped space, $X_k$ to be  
$X_k = \Gamma(G,{\mc A}) \cup (\cup \{{\mc H}(\Gamma_{i,t})_k \mid 1\leq i\leq n,
t\in T_i\})$.

\subsection{Relatively hyperbolic groups and first algorithms}

The cusped space can be used to define {\em relative hyperbolicity}:

\begin{defi} \label{d:RH}
Suppose that $G$ is a finitely generated group and that 
${\bf P} = \{ P_1, \ldots , P_n \}$ is a finite collection of finitely
generated subgroups of $G$.  Suppose further that $\mathcal A$
is a finite compatible generating set for $G$.

Then $G$ is hyperbolic relative to ${\bf P}$ if the cusped space associated
to $(G,{\bf P},\mathcal A)$ is $\delta$-hyperbolic for some $\delta$.
\end{defi}

It is proved in \cite[Theorem 3.25]{GrovesManning} that this definition
is equivalent to other standard definitions. 

\begin{assumption} \label{a:properlyRH}
Since the main results of this paper are straightforward in case
$G$ is equal to a parabolic subgroup, we 
will always assume that each parabolic subgroup of $G$ is properly contained
in $G$.  Since a finite subgroup can be omitted from the list of parabolics without
affecting relative hyperbolicity, we will also assume that parabolic subgroups are 
infinite.
\end{assumption}

\begin{definition} \label{def:Boundary}
Let $G$ be a relatively hyperbolic group, with associated cusped space $X$.  The
boundary of $G$, denoted $\partial G$, is the Gromov boundary of the space
$X$ (defined as usual by equivalence classes of quasi-geodesic rays).
\end{definition}

It follows from \cite[Sections 6,9]{Bowditch_rel_hyp} that the above definition
of $\partial G$ agrees with Bowditch's definition (the analysis in 
\cite[Section 6]{Bowditch_rel_hyp} applies to our cusped space $X$).  In
particular, the results of Bowditch quoted in Section \ref{s:Splittings}
below apply to the boundary defined in Definition \ref{def:Boundary}.

 Let us now recall the main result of \cite{Dfind}.

\begin{theorem} \cite[Theorem 0.2 and 0.1]{Dfind} \label{t;dfind}

There is an algorithm that, given a finite presentation of a group $G$, terminates if and only if $G$  is hyperbolic relative to abelian subgroups.  In case it terminates, it provides a finite presentation of each of the parabolic subgroups (up to conjugacy) in terms of the given generators of $G$, together with a constant for a linear relative isoperimetric inequality (or equivalently  a linear combinatorial isoperimetric
inequality of the coned-off Cayley complex, see
\cite[Section 2]{GrovesManning}).
\end{theorem}

\begin{remark}
Given a finite presentation of an abelian group, it is straightforward to determine
if the group is finite or not.  Therefore, it is no trouble to make the Assumption
\ref{a:properlyRH} that all parabolic subgroups are infinite.  

If a parabolic subgroup is virtually cyclic, we can exclude it from the list of parabolics
without affecting relative hyperbolicity.  Since we can determine if an abelian group
is virtually cyclic from a finite presentation, we can also assume that all parabolic
subgroups are not virtually cyclic.
\end{remark}

We require one further result which will allow us 
to effectively compute a value of $\delta$ for which the cusped space is
$\delta$-hyperbolic.  

\begin{theorem}\cite[Theorem 3.24]{GrovesManning} \label{t:ConetoCusp}
Suppose that $G$, ${\bf P}$ and $\mathcal A$ are as in Definition \ref{d:RH},
and that $G$ is hyperbolic relative to ${\bf P}$.  If the coned-off Cayley complex
of $G$ associated to $G, {\bf P}, {\mc A}$ satisfies a linear combinatorial isoperimetric
inequality with constant $K$, then the cusped space associated to $G, {\bf P}, {\mc A}$
satisfies a linear isoperimetric inequality with constant $K_1 = 3K(2K+1)$.
\end{theorem}

The space $X$ is the cusped space from Definition \ref{d:X1}. We can compute $\delta$
using  Theorems  \ref{t;dfind} and  \ref{t:ConetoCusp}.

Let $v_0$ be the natural base point of $X$, corresponding to the identity element of $G$ which
is a vertex of the Cayley graph of $G$, which is embedded in the cusped space $X$. 

 Balls of $X$ can be computed, by using a solution to the word problem in $G$
 (and the fact that we know we have a compatible generating set, and know which
 of our generators lie in which parabolic subgroups -- 
 this uses the fact that the parabolics are abelian).

\begin{lemma} \label{lem:geod_ray}
For every point $x \in X$ there is a geodesic
ray which starts at $v_0$ and passes within $3\delta$ of $x$.
\end{lemma}
\begin{proof}
Suppose that $x \in X$ has the form $(\gamma,k)$, where $\gamma \in t_j P_i$,
for some parabolic subgroup $P_i$ and some transversal element $t_i$.  (Note that
in case $k = 0$, there may be more than one such pair $(P_i,t_j)$.  This is unimportant
in the following argument.)

Let $\alpha \in \partial G$ be the point in the boundary fixed by $t_j P_i t_j^{-1}$,
corresponding to the horoball built from $t_j P_i$.  Since $P_i$ is a proper subgroup
of $G$, there exist $\beta \in \partial G$ so that $\beta$ is the fixed point of
some $t_k P_i t_k^{-1}$, and $t_k P_i \neq t_j P_i$.

Let $p$ be the path which consists of a shortest path between the horoballs
$t_j P_i$ and $t_k P_i$, and two vertical paths to $\alpha$ and $\beta$.  It is
easy to see that $p$ is a bi-infinite geodesic path.
By applying an element of $t_j P_i t_j^{-1}$ to the path $p$, we may assume that
$p$ passes through $x$.

According to \cite[Lemma 2.11]{GrovesManning} (which is a simple exercise), a geodesic
triangle in a $\delta$-hyperbolic space with some or all vertices ideal is
$3\delta$-slim.  Consider a partially ideal triangle with vertices $v_0,\alpha,\beta$,
where the path between $\alpha$ and $\beta$ is $p$.  Then $x$ lies within at
most $3\delta$ from one of the edges $[v_0,\alpha]$, $[v_0,\beta]$, as required.
\end{proof}

\subsection{Constants} \label{ss:Constants}

We choose $\delta$ to be a positive integer so that $X$ is $\delta$-hyperbolic.

Let $C = 3\delta$, so by Lemma \ref{lem:geod_ray} any point in $X$ is at distance at most $C$ from a
geodesic ray starting at $v_0$.

Let us define $M = 6(C+45\delta) + 2\delta +3$, 
and $k = 2 M$, and recall that $X_k$ is the $k$-thick part of the cusped space. 
Let $K = 3  (2^{2M+3}) + M + 3$.
Finally, let $R(n) = 4(n+M) +3k +50 \delta +3$.

\subsection{Two properties of horoballs} \label{ss:HoroProps}

Recall that in any metric space $(X,d)$ the {\em Gromov product} of points
$x$ and $y$ with respect to a third point $z$ is
\[      (x \cdot y)_z = \frac{1}{2} \left( d(x,z) + d(y,z) - d(x,y) \right)     .       \]
Given a comparison tripod $Y_{xyz}$ for $x,y,z$ (where the vertices
of the tripod are $\bar{x}, \bar{y}, \bar{z}$), the distance from $\bar{z}$ to the
centre of the tripod is $(x \cdot y)_z$.  In particular, if $z$ lies on a geodesic between
$x$ and $y$ then $(x \cdot y)_z = 0$.

A definition of $\delta$-hyperbolicity is that for all $w,x,y,z \in X$ we have
\[      (x \cdot y)_w \ge \min \{ (x \cdot z)_w, (y \cdot z)_w \} - \delta      .       \]
(Because the constants for translating between the various definitions of $\delta$-hyperbolicity are explicit, we assume that the above equation holds as well as 
triangles being $\delta$-thin, and $\delta$-slim.  See \cite[III.H]{BH} for more details.)

 \begin{lemma}\label{lem;horoballs1}

        Let $H_k$ be an 
        horoball at depth $k>0$,  $p_0\in H_k$, and 
        $p_1\in H_{k-1}\setminus H_k$. 
        
        Then, there exists $p_2 \in  H_{k-1}\setminus H_k$  such that 
        $d(p_0,p_2)\leq d(p_0,p_1) + 3$, 
        and $(p_1 \cdot p_2 )_{p_0} \le  3$.

\end{lemma}
\begin{proof}
Suppose that the horoball is built on the graph $\Gamma$, so its vertices
are labelled by $(\gamma,i)$ where $\gamma \in \Gamma^{(0)}$ and $i \ge k-1$.
Note that we are assuming that $\Gamma$ has infinite diameter.

We consider the geodesic $\gamma$ between $p_0$ and $p_1$.  By Proposition 
\ref{p:HoroProps}.\eqref{item:geo}, we may assume that $\gamma$ 
either consists of a single
vertical segment, or else at most two vertical segments and a single horizontal
segment of length at most $3$.
We distinguish three cases:  where $\gamma$ contains two vertical segments (and one horizontal segment); where
$\gamma$ is entirely vertical; and where
$\gamma$ contains a single vertical segment and a horizontal segment.

Suppose first that this geodesic contains a horizontal segment and two vertical segments, and suppose that
$p_0 = (\gamma,i)$.  Then the path formed by concatenating the geodesic from
$p_1$ to $p_0$ and the vertical path from $p_0$ to $(\gamma,k-1)$ is a geodesic,
and $(p_1 \cdot (\gamma,i))_{p_0} \le \frac{1}{2}\delta$, so we may take $p_2 = (\gamma,k-1)$.
In this case certainly we have $d(p_0,p_2) \le d(p_0,p_1)$.
Note also that since $p_0$ lies on a geodesic from $p_1$ to $p_2$ we must
have $(p_1 \cdot p_2)_{p_0} = 0$.

Next suppose that the geodesic from $p_1$ to $p_0$ is vertical, so
$p_1 = (\gamma,k-1)$ and $p_0 = (\gamma, i)$ for some $\gamma \in \Gamma^{(0)}$
and some $i > k-1$.  Note that if $d(p_0,p_1) \le \delta$ then we can take $p_2 = p_1$
and the lemma is trivial.  Thus we may suppose that $i - k-1 > \delta$.

Since the graph $\Gamma$ has infinite diameter, there is a vertex 
$\gamma_0 \in \Gamma^{(0)}$ so that $d_{\Gamma}(\gamma,\gamma_0)
= 3(2^{i})$.  Thus the distance between $(\gamma,i)$ and $(\gamma_0,i)$
at depth $i$ is $3$, and there is a horizontal path between them which is
a geodesic in $H_k$.  The path between $(\gamma,k-1)$ and
$(\gamma_0,k-1)$ which consists of the vertical paths between $(\gamma,k-1)$
and $(\gamma,i)$ and between $(\gamma_0,k-1)$ and $(\gamma,k-1)$ and
the horizontal path between $(\gamma,i)$ and $(\gamma_0,i)$ is a geodesic
in $H_k$.  Let $p_2 = (\gamma,k-1)$.  Clearly
$d(p_0,p_2) = d(p_0,p_1) + 3$, and it is not hard to see that
$(p_1 \cdot p_2)_{p_0} = 0$ (since again $p_0$ lies on a geodesic between $p_1$ and $p_2$).

Finally, suppose that the geodesic from $p_1$ to $p_0$ contains a single vertical segment and a horizontal segment.  In this case the geometry of the graph on which
the horoball is based may prevent us from finding a point $p_2 \in H_{k-1} \setminus
H_k$ so that $p_0$ lies on a geodesic from $p_1$ to $p_2$.  However,  let
$x$ be the point at the end of the vertical segment in the geodesic from $p_1$ to $p_0$.
Then $d(x,p_0) \le 3$ and $d(p_0,p_1) = d(p_0,x) + d(x,p_1)$.  
The second case above finds a point $p_2$ so that $x$ 
lies on the geodesic from
$p_1$ to $p_2$ and $d(x,p_2) = d(x,p_1) + 3$.  Since $x$ lies on the geodesic from
$p_1$ to $p_2$ a simple calculation then shows that 
$(p_1 \cdot p_2)_{p_0} \le d(x,p_0) \le 3$.
Also, $d(p_0,p_2) \le d(p_0,x) + d(x,p_2) = d(p_0,x) + d(x,p_1) + 3 = d(p_0,p_1) + 3$,
as required.

This finishes the proof in all cases.
\end{proof}

\begin{definition}
Suppose that $P$ is a finitely generated abelian group.  A finite generating 
set $A$ for $P$ is called {\em sensible} if it can be partitioned into two subsets
$A = A_1 \sqcup A_2$ where $A_1$ is a basis for a free abelian group and $A_2$ is a generating set for a finite abelian group $P_f$, and $|A_2|$ is minimal amongst all
generating sets of $P_f$.
\end{definition}

\begin{remark}
The algorithm from \cite{Dfind} which finds the relatively hyperbolic structure for a
relatively hyperbolic group with abelian parabolics finds a sensible generating
set for the parabolics.  Therefore, whenever $G$ is assumed to have abelian parabolics, we will assume that the generating set $\mc A$ is such that the intersection of 
$\mc A$ with any element of ${\bf P}$ is sensible.
\end{remark}

\begin{lemma} \label{lem:path}
Suppose that $A$ is a sensible generating set for an infinite abelian group $P$ which is not virtually cyclic.  For any $r \ge 1$, if $a,b \in P$ are at distance at least $r$ from
 $1$  then there is a path of length at most $3d(a,b)$ from
 $a$ to $b$ which does not intersect the ball of radius $r-1$ about $1$ (all distances are measured with respect to the word metric coming from the sensible generating set).
\end{lemma}
\begin{proof}
Let $A = A_1 \sqcup A_2$ where $A_1$ generates a free abelian subgroup $P_1$
of $P$, and $A_2$ generates a finite group.  
Since $P$ is not virtually cyclic,  $P_1$ has rank at least $2$.

If $d(a,b)\leq 1$, there is nothing to prove.

Let us write $a=(a_1,a_2)$, with $a_i  \in \langle A_i\rangle$.
   and $b=(b_1,b_2)$ similarily. Let us consider the coordinates of $a_1=(x_1,\dots, x_m)$ and $b_1=(y_1,\dots, y_m)$ in the basis $A_1$ (we choose additive notations). If two coordinates corresponding to a basis element $e\in A_1$  are of the same sign (assume it is positive, up to changing $e$ into $-e$), then   $e+a$ and $e+b$ are both at distance at least $r+1$ from $1$. If  $\sigma$ be a geodesic segment between $a$ and $b$, $e+\sigma$ has its extremal edges outside the ball of radius $r$ about $1$. Considering the neighbors of $e+a$ and $e+b$ on    $e+\sigma$, we get two points $a'$ and $b'$ satisfying the assumption of the lemma, with $d(a',b') = d(a,b) -2$, and it takes two paths of length $2$ to reach them from $a$ and $b$ outside  the ball of radius $r-1$ about $1$.

Thus by iterating this process, we can reduce the case to that where all coordinates of $a_1=(x_1,\dots, x_m)$ have opposite signs than that of $b_1=(y_1,\dots, y_m)$.

Let us write   $\langle A_2\rangle$ as a direct product of cyclic groups, and $a_2 = (\alpha_1,\dots  ,\alpha_k)$ and $b_2 = (\beta_1,\dots  ,\beta_k)$ in this product. If, in the cyclic component, $0$ is not on a geodesic between $\alpha_i$ and $\beta_i$, and assuming that $a$ is closer to $1$ than $b$, then there is a path of length $\|\alpha_i - \beta_i\|$ between $a$ and $a'$, not approaching $1$, such that the coordinate of $a'$ along this cyclic factor is that of $b$.  The distance to $b$ has reduced by the same amount, thus, we can reduce to the case where, in  $\langle A_2\rangle$, a geodesic between $a_2$ and $b_2$ goes through $0$, that is $ \| a_2-b_2 \| = \| a_2\| + \|b_2 \|  $.

Up to changing some elements $e$ of the basis  $A_1$ into $-e$, we can assume that all coordinates of $a_1$ are negative, and all that of $b_1$ are positive. Then, by adding $(0,\dots,0,-1) \in \langle A_1\rangle  $, then  $(1,0,0,,\dots,0)  \in \langle A_1\rangle$ to $a$, exactly $|x_1| $ times, one defines a path of length $2 |x_1|$, outside the ball of radius $r-1$, to a point with $0$ first coordinate. 
 Thus, we can find a path from $a$ to $( (0, \dots, -\sum_1^{m} | x_i|  -  \| a_2\|), 0)$ of length at most $ 2 \sum_1^{(m-1)} | x_i|$, outside the ball of radius $(r-1)$. Similarily, one can find a path from $b$ to $ ( (\sum_1^{m} y_i + \| b_2\|   ,0, \dots, 0),0)  $, of length at most  $ 2 \sum_2^{m}  y_i + 2\|b_2\|$,  outside the forbidden ball. These two paths have total length at most  $ 2\sum_1^{m}  (y_i-x_i) + 2\| a_2 -b_2\| $ (recall $x_i<0$, and $ \| a_2-b_2 \| = \| a_2\| + \|b_2 \|  $ ) , thus at most $2d(a,b)$.
 Notice now that there is a path between $( (0, \dots, -\sum_1^{m} | x_i| - \| a_2\|), 0)$ and  $ ( (\sum_1^{m} y_i +\| b_2\|,0, \dots, 0),0)$, via  $ ( (\sum_1^{m} y_i +\| b_2\| ,0, \dots, 0, -\sum_1^{m} | x_i|- \| a_2\|  ),0)$, and of length  $\sum_1^{m}  y_i-x_i + \| a_2\| + \|b_2 \|  = d(a,b)$. This path does not get closer to $1$ than its end points.  To summarise, there is a path of length at most $3d(a,b)$ from $a$ to $b$ outside the forbidden ball. 
\end{proof}

Recall that we defined $K = 3 (2^{2M+3}) + M + 3$.

 \begin{lemma}\label{lem;horoballs2}

   Assume that the parabolic subgroups of $G$ are abelian, and that the intersection
   of $\mc A$ with each element of ${\bf P}$ forms a sensible generating set.

        Let $H_k$ be an 
         horoball at depth $k>0$, and $x,y$ in $H_k$ be such that $d(x,y) \leq M$, and  
         $|d(x,v_0) -d(y,v_0)| \leq 20\delta$.

         Then, there exists  a path between $x$ and $y$ in $H_k$ which is of length
        at most $K$ which does not intersect the ball
        of radius $ m =\min \{ d(v_0,x) , d(v_0,y) \}$ about $v_0$.
\end{lemma}
\begin{proof}
Let $\sigma_1$ be a geodesic from $v_0$ to $x$ whose intersection with $H_k$
satisfies the conclusions of
Proposition \ref{p:HoroProps}.\eqref{item:geo}, and $\sigma_2$ such a geodesic
from $v_0$ to $y$.
Since, for all $i$,  $H_i$ is  convex in $X$,  the points at which $\sigma_1$ and $\sigma_2$
first intersect $H_1$ are distance at most $2\delta$ apart, and since $k>
1+\log_2(2\delta)$, the 
points at which $\sigma_1$ and $\sigma_2$
first intersect $H_{k}$ are distance at most $1$ apart.

The paths $\sigma_1$ and $\sigma_2$ then travel vertically at distance at most $1$
apart, then one of
the paths turns horizontal, and goes back up, whilst the other possibly goes
deeper, before possibly turning and going up also.  If both paths contain only a
single vertical segment,
the lemma is straightforward:  suppose $x$ is the higher of $x$ and $y$.  Then go 
directly downwards from $x$ to the depth of $y$, and then apply Lemma \ref{lem:path}.
Thus we can assume that $\sigma_1$ contains two vertical segments.

Suppose (by relabelling $x$ and $y$ if necessary) that $x$ lies no deeper than $y$, and
that the depth of $x$ is $i$.
Since $x$ and $y$ lie no more than $M$ apart in $X$, we may append a path 
$\gamma$ to the end of $\sigma_2$ so that $\sigma_2 \cup \gamma$ is a geodesic 
and so that $\gamma$:
\begin{enumerate}
\item lies outside of the ball of radius $m$ about $v_0$;
\item is entirely vertical except possibly for a horizontal segment of length at most $3$;
\item has length at most $M + 3$; and
\item has end $y'$ (at the other end of $y$) which is at depth $i$ (the same as $x$) and
is at distance at most $2M+3$ from $x$ in $X$.
\end{enumerate}

We will find a path from $y'$ to $x$ of bounded length which stays outside of the
ball of radius $m$ about $v_0$.  Let $\sigma_2'$ be the appended geodesic  $\sigma_2 \cup \gamma$.

Let $d_i$ denote distance at depth $i$.
Since the distance in $X$ between $x$ and $y'$ is at most $2M+3$, we have
$d_i(x,y') \le 2^{2M+3}$.

The intersection of level $i$ with $\sigma_1$, consists of two
points, $x$ and another point $p_1$.  The intersection of $\sigma_2'$ with level
$i$ consists of points $y'$ and $p_2$.  Note that, since $i>k$, as noted before, one has $d(p_1,p_2) \le
1$, hence  $d_i(p_1,p_2) \le 1$.

Now, the part of $H_k$ at level $i$ is a Cayley graph of a finitely generated
abelian group $P$, where $P \in {\bf P}$.  This Cayley graph is with respect
to a  generating set consisting of all elements of length at most $2^i$ in a 
sensible  generating set of $P$.  Let $L_i$ be this Cayley graph at depth $i$ in $H_k$.
 By Lemma \ref{lem:path}, for any points $a$ and $b$ in $L_1$ which lie
at distance at least $r$ from $1$, there is a path of length at most $3 d_1(a,b)$
between $a$ and $b$ which does not intersect the ball of radius $r-1$ about $1$.
Since the distance in $L_i$ between points $u = (u_0,i) ,v = (v_0,i) \in L_i$ is 
$\lfloor 2^{-i}d_1(u_0,v_0) \rfloor$,  it is easy to see that Lemma \ref{lem:path}
holds for $L_i$ also.

Let $\sigma$ be a path in $L_i$ between $x$ and $y'$ which does not come closer
to $p_1$ and $p_2$ than $x$ or $y'$, and is of minimal length with respect to
this restriction. By the remark above, the length of $\sigma$ is at most   $3 d_i(x,y')\leq  3 ( 2^{2M+3}) $.

In conclusion, there is a path between $x$ and $y'$ of the right kind of length
at most $3 ( 2^{2M+3})$, and the path between $y$ and $y'$ has length at
most $M+3$.  This finishes the proof of the lemma.
\end{proof}

\begin{remark}
There are versions of Lemma \ref{lem;horoballs2} with similar conclusions and
much weaker assumptions on the parabolic subgroups.  However, we were unable
to prove a result with no assumptions on the parabolics, and different types of paths of different lengths and with slightly different conclusions are needed for different types of parabolics.  Thus, we decided to present 
a simple version which suffices for the needs of this paper.
\end{remark}

\section{Splittings} \label{s:Splittings}

        In this section we gather some results due to B.~Bowditch about splittings.
         We first need a characterisation of connectedness of the boundary.

        \begin{proposition} \cite[Propositions 10.1-3]{Bowditch_rel_hyp} \label{propB;one-end=connected}
        
         Let $(G,{\bf P})$ be a relatively hyperbolic group. 
          Its boundary $\partial G$ is disconnected if and only 
          if $G$ splits nontrivially over a finite group relative to
          ${\bf P}$. In this case, every vertex group is  hyperbolic
          relative to the parabolic subgroups it contains.

        In particular, if the parabolic subgroups of $G$ are all one-ended, then 
        $G$ is one-ended if and only if $\partial G$ is connected.

        \end{proposition}
        In the above result, a splitting is {\em relative to ${\bf P}$} if each element
        of ${\bf P}$ is elliptic in the splitting.

        Then we need to characterise the presence of global cut points. 
        For that, we recall the notion of peripheral splitting.

        \begin{defi} \cite{Bowditch_periph}

        Let $(G,{\bf P})$ be a relatively hyperbolic group. 
        A {\em peripheral splitting} of $G$ is a representation of $G$ as 
        a finite bipartite graph of groups where ${\bf P}$ consists precisely 
        of the (conjugacy classes of) vertex groups of one colour.  
        A peripheral splitting is a {\em refinement} of another if 
        there is a colour preserving 
        folding of the first into the second. 
        \end{defi}

        \begin{theorem} (B. Bowditch) \label{theo;bowditch}

          Let $(G,{\bf P})$ be a relatively hyperbolic group. Assume that $\partial G$ is connected.

          \begin{enumerate}

          \item \cite[Theorem 0.2]{Bowditch_connected}   If every maximal parabolic subgroup of $G$  
                is (one or two)-ended, 
                finitely presented, and  without infinite torsion subgroup, then, every global cut point 
                of $\partial G$ is a parabolic fixed point. 

          \item \cite[Theorem 1.2]{Bowditch_boundGF} If there is a global cut point that is a parabolic point, 
            then $G$ admits a proper peripheral splitting.

          \item \cite[Theorem 1.2]{Bowditch_periph} If $G$ admits a proper peripheral splitting, 
                then $\partial G$ admits a global cut point.

          \end{enumerate}
                  
        \end{theorem}

        It is established in \cite[Theorem 1.3]{Bowditch_periph} that the vertex groups of peripheral splittings are hyperbolic relative 
        to the trace of the parabolic subgroups of $G$ on them, and that  if $\partial G$ is connected, then, 
        any non-peripheral vertex group also has connected boundary, and is hyperbolic relative to its adjacent edge groups.

        We also need an accessibility result.

        \begin{theorem}\cite[Theorem 1.4]{Bowditch_periph}. \label{theoB;access}

          Suppose that $G$ is relatively hyperbolic with connected boundary.
          Then $G$ admits a (possibly trivial) peripheral splitting which is maximal in the sense
          that it is not a refinement of any other peripheral splitting.  
  
          The boundaries of the components of such a maximal splitting 
          do not contain any global cut point.  

        \end{theorem}

We will recognise when the boundary is connected without global cut point, on the one hand, and when the group admits nontrivial splittings on the other.

\section{On a condition of M. Bestvina and G. Mess} \label{s:BM}

Let $G$ be a relatively hyperbolic group, and $X$ the cusped space defined
in Section \ref{s:Prelim}.

Let us introduce a slight variation of the property $\ddagger$ of
Bestvina and Mess, defined  in \cite{BestvinaMess}.

Given  $\epsilon \geq 0$, and two points $x,y$ in $X$, we say that $x$ and $y$ satisfy $\star_\epsilon$ if 
\[
   \star_\epsilon : \qquad   | d(v_0,x) - d(v_0,y)| \leq \epsilon   \hbox{ and } d(x,y)\leq M
\]

Given an integer $n$, we say that two points $x,y$ satisfying $\star_\epsilon$, satisfy $\ddagger(\epsilon,n)(x,y)$ if there is a path of length at most $n$ from $x$ to $y$, in the complement in $X$ of the ball $B_{v_0}(m-C-45\delta + 3 \epsilon)$, where
$m = \min \{ d(v_0,x), d(v_0,y) \}$.

First one can reproduce {\it verbatim} the argument of Bestvina and Mess (we refer the reader to \cite[Proposition 3.2]{BestvinaMess}) to get the following property.

 \begin{lemma}\label{lem;BM-loc-connected}
    If there exists  $n$ 
        such that  for all $x,y$ 
    in $X$ satisfying $\star_0$,  the property $\ddagger(0,n)(x,y)$ holds, then
    $\partial G$ is connected.
  \end{lemma}

Note that if  $\ddagger(0,n)(x,y)$ holds, then $\ddagger(0,m)(x,y)$ 
holds for all $m\geq n$ also.
Following another argument of Bestvina and Mess (\cite[Proposition 3.3]{BestvinaMess}), we also collect (recall that we defined the constant $k$ in Subsection \ref{ss:Constants}):

\begin{lemma} \label{lem;BM-cutpoint}

 If the boundary $\partial G$ is connected, 
and has no global cut point, 
then there is    $n$ 
 such that $\ddagger(10\delta,n)(x,y)$ holds for all $x,y$ in $X_k$.
\end{lemma}
\begin{proof}
We prove the contrapositive. If for each $n$ there are $x_n, y_n$ so that
$\ddagger(10\delta,n)(x_n,y_n)$ is false, we will find a 
 point $\xi$ in $\partial G$ such that $\partial G
\setminus \{\xi\}$ is disconnected.

Thus, let us  assume that $\ddagger(10\delta,n)(x_n,y_n)$ is false for a sequence of
vertices $x_n, y_n$ in $X_k$. Let $r_n$ and $s_n$ be
geodesic rays from $v_0$ that pass at distance at most $C$ from $x_n$ and
$y_n$, respectively. Their points at infinity are denoted by $\alpha_n$ and $\beta_n$.

Since there are finitely many orbits of vertices in $X_k$, one can find a
sequence of elements $\gamma_n$ such that, after extraction of a subsequence,
$\gamma_n x_n$ is constant, say equal to $x$. Since $y_n$ is at bounded
distance of $x_n$, and the space $X$ is proper, we can assume that $\gamma_n
y_n$ is constant, say equal to $y$, with $d(x,y) \leq M$. Without loss of
generality, one can assume that $\min\{d(v_0,x_n),d(v_0,y_n) \} =
d(v_0,x_n)$, for all $n$ in the extracted subsequence.  We also extract in
order that $\gamma_n \alpha_n \to \alpha$, and $\gamma_n \beta_n \to \beta$ in 
 $\partial X$, $\gamma_n v_0\to \xi \in \partial G$, $[x,\gamma_n v_0] \to
[x,\xi)= \rho$   and
$r_n \to r$, and $s_n \to s$, bi-infinite geodesics from $\xi$ to $\alpha$,
and $\beta$ (it is easy to check that $\xi$ is different from $\alpha$ and
$\beta$).

We consider then  $B_n = \gamma_n B(v_0, d(v_0,x_n) -
15\delta -C)$.

 Let $B = \bigcup_{\mathbb{N}} B_n$. Clearly $\xi \in \bar{B}
 \cap \partial X$ for the usual topology on $X \cup \partial X$. One has furthermore $\{\xi\} =  \bar{B}
 \cap \partial X$. For,  if $\zeta \in   \partial X \setminus \{ \xi \}$ then, letting $\{ f_n \}$
be a sequence of points in $B_n$ going to $\zeta$, by hyperbolicity  in the pentagons $(x,\xi,\gamma_n v_0, f_n, \zeta)$,  there is $N>0$ such that if $n$ is large enough, any path from $\gamma_n v_0$ to $f_n$ passes at distance at most $N$ from $x$. This provides $d(  \gamma_n v_0, f_n) > 2 (d(v_0,x_n) -
15\delta -C )  $, for $n$ large enough, contradicting the assumption that both are in $B_n$.

 The space $X\setminus B$ is disconnected. Indeed, $x$ and $y$ are in $X\setminus B$, and 
  if there is a path $\lambda$ in $X \setminus B$, from $x$ to
  $y$, then, for $n>length(\lambda)$, translating by $\gamma_n^{-1}$, we see that   $\ddagger(10\delta,n)$ is true for $(x_n,y_n)$, contradicting the assumption on $x_n, y_n$.

  The distance between $x$ and $r$ is at most $C$, and that between $x$ and
  $B$ is  $15\delta +C$. Therefore the ray $r$ enters the component of $x$ in
  $X\setminus B$. 
  Similarly, the ray $s$ enters the component of
  $X\setminus B$ of $y$.

  For a component $A$ of $X\setminus B$,   
  let us consider $\mathcal{O}(A) \subset
  \partial X \setminus\{\xi\}$, the set of the points $\zeta$ such that a 
  bi-infinite
  geodesic $(\xi, \zeta)$ enters $A$.

  The set $\mathcal{O}(A)$ is
  open: indeed, if $\zeta$ is in it, a neighborhood of $\zeta$ consists of
  points $\eta$ such that $(\xi,\eta)$ remains at distance $10\delta$ from
  $(\xi,\zeta)$, from $\xi$ until some point at distance at least $20\delta$
  from $B$, and therefore these geodesics must enter 
  the same component of $X\setminus B$ as $(\xi,\zeta)$.

  If $A' \neq A$, $\mathcal{O}(A)$ and $\mathcal{O}(A')$ are disjoint, because
  if a point $\zeta$ is in the intersection, one geodesic $(\xi,\zeta)$ enters
  $A$, although another enters $A'$, but both remain at distance $2\delta$,
  thus providing paths between $A$ and $A'$ that avoid 
  $B$ (contradicting that $A\neq A'$).

  Finally there are at least two 
  non-empty such open sets (one containing $\alpha$, and 
  one containing $\beta$), and the union of them covers 
  $\partial X \setminus \{\xi\}$. 
  All this proves that $\partial G \setminus\{\xi\}$
  is disconnected. 
\end{proof}

We now want to reduce  $\ddagger(0,n)$ on $X$,  to $\ddagger(10\delta,n)$ on some finite set. We begin with two lemmata, dealing with the $k$-thick part $X_k$.
Recall the definition of $R(n)$ from Subsection \ref{ss:Constants}.
 
\begin{lemma}\label{lem;good_gamma}

        For all $x$ in $X_k$ with 
        $d(x, v_0) \geq R(n)$, 
        if $p$  is a point on a geodesic $[v_0, x]$ at distance $2 (n+M)$  from $x$, 
        then there exists an element $\gamma$ such that   
        $(\gamma v_0 \cdot x)_{p}\leq \delta$, and 
        $ d(x,\gamma v_0)    \leq R(n) -50\delta$. 

\end{lemma}

  Notice that this is straightforward for hyperbolic groups, since in the Cayley graph, one can
always choose $\gamma v_0=p$.

\begin{proof}
Let us  consider the subsegment $[p',p]$ of $ [v_0, p]$ and of length $4k$  
 (it exists since $d(v_0,x) \geq R(n)$).
If the distance between $p'$ and $G v_0$ is at most $3k$, then, the element $\gamma$ such that $\gamma v_0$ is closest to $p'$ is suitable.

  If now $p'$ is at distance at least $3k$ from $G v_0$, then it is not in $X_k$. 
Since $x\in X_k$, there is a first point  $p''$  on $[p', x]$ that is in $X_k$.

By assumption on $p'$, one has $d(p',p'')\geq 2k$.
 Moreover since the horoballs at depth $k$ are separated, the segment $[p',p'']$ remains in an horoball. By Lemma \ref{lem;horoballs1}, there is another point $q$ in $X_k$ with $d(p',q) \leq d(p',p'') +3$ and such that 
$(q \cdot p'')_{p'} \leq  3$.

Since $q \in X_k$, there exists $\gamma \in G$ such that $d(\gamma v_0, q ) \leq k$. Recall that   
        $d(p', \gamma v_0) \geq 3k$  hence $d(p',q) \geq 2k$. Therefore, 
for this element, $(\gamma v_0 \cdot q)_{p'} \geq  (d(p',q) + d(p', \gamma v_0)  - k)/2 \geq  2k$, 
and by hyperbolicity, it follows that $(\gamma v_0 \cdot p'')_{p'} \leq \delta + 3$. Since $(p \cdot p'')_{p'} = d(p',p'') \ge 2k$, one has   $(\gamma v_0 \cdot p)_{p'} \leq 2\delta + 3$.  
Therefore,  $(\gamma v_0 \cdot p')_{p} \geq d(p,p') -  (2\delta + 3) >\delta$. Since $(x \cdot p')_{p} =0$, by hyperbolicity, $(\gamma v_0 \cdot x)_{p}\leq \delta$.

Moreover, one has $d(x,\gamma v_0) \leq d(x,p')+d(p',q)+d(q, \gamma v_0)$.
Recall that  $d(p',q)\leq d(p',p'')+3\leq d(p',x)+3$, and that $ d(q, \gamma v_0) \leq   k$. Since  $d(p',x)\leq d(p,x) +k$, 
one has $ d(x,\gamma v_0)    \leq  2d(x,p)+6k +3$ that is, at most $4(n+M) +3k +3  \leq R(n) - 50\delta$.  Therefore, $\gamma v_0$ is suitable.
\end{proof}

 \begin{lemma} \label{lem;move}

 Assume that,  for some  $n$, $\ddagger(10\delta,n)(x,y)$ holds for all pairs of vertices $(x,y)$ in $ B_{v_0} (  R(n)  )$ satisfying   $\star_{10\delta}$.

   Let $x,y \in X_k$, satisfying $\star_0$,  and $d(v_0,x) \geq R(n)$. 
        Let $p$ on a geodesic segment  $[v_0,x]$ at distance 
        $2(n  +M)$    from $x$, as  in the previous lemma. 
  Then the pair $(x,y)$ satisfies $\ddagger(0,n)$.  
 \end{lemma}
\begin{proof}
By Lemma \ref{lem;good_gamma}, there exists $\gamma \in G$ such  that   
   $(\gamma v_0 \cdot x)_{p} \leq    \delta$ and $d(\gamma v_0, x) \leq R(n)-50\delta$. 
Let us first prove that $(\gamma^{-1}x , \gamma^{-1} y)$ satisfies $\star_{10\delta}$.

By hyperbolicity, $\min\{(\gamma v_0 \cdot  y)_{p} , (x \cdot  y)_{p}  \} \leq (\gamma v_0 \cdot  x)_{p}+\delta \leq 2\delta$.
 Since $(x \cdot  y)_{p} \geq M+2n$, it follows $(\gamma v_0 \cdot  y)_{p}\leq 2\delta$, and therefore, 
$d(\gamma v_0, y) \geq d(\gamma v_0, p) + d(p, y) -6\delta.$
Since $(\gamma v_0 \cdot  x)_{p}\leq \delta$, similarly, $d(\gamma v_0, x) \geq d(\gamma v_0, p) + d(p, x) -2\delta$.  

By triangular inequality, on the other hand, $d(\gamma v_0, y) \leq d(\gamma v_0, p) + d(p, y) $, and $d(\gamma v_0, x) \leq d(\gamma v_0, p) + d(p, x)$
One can thus estimate $ |d(p, x) - d(p, y)| - 2\delta \leq |d(\gamma v_0,x) - d(\gamma v_0, y)| \leq |d(p, x) - d(p, y)| + 6\delta$.

  By assumption $\star_0$ on $x$ and $y$, and since  $p$ is on $[v_0,x]$ and  $\delta$-close of  $[v_0,y]$ , one has $|d(p, x)-d(p, y)| \leq 2\delta$. 
This ensures that $(\gamma^{-1}x , \gamma^{-1} y)$ satisfy $\star_{10\delta}$.

  By assumption, $d(\gamma v_0,x) \leq  R(n)-50\delta$, hence 
there exists a path $c: [0,n] \to X$ from $x$ to $y$, such that for all $t$,
$d( \gamma v_0, c(t)) \geq  \min \{d(\gamma v_0, x),d(\gamma v_0, y)  \} -C -
15\delta$, in particular,  $d (\gamma v_0, c(t))  \geq  d(\gamma v_0, x) - C-25\delta$.

We  now need to control $d(v_0,c(t))$ for all $t$.
Changing $y$ by $c(t)$   in the argument above provides  the estimation  
$    d(p, c(t) ) - d(p, x) \geq  d(\gamma v_0, c(t))- d(\gamma v_0,x) - 6\delta$. 
 
One also easily gets that $d(v_0,c(t)) - d(v_0,x) \geq    d(p, c(t))-d(p, x)   -2\delta$, thus, 
 $d(v_0,c(t)) - d(v_0,x) \geq   d(\gamma v_0, c(t)) - d(\gamma v_0,x)
    -6\delta\geq - C-35\delta$, as required.
    \end{proof}

\begin{corollary}\label{coro;allpairs}

        Let $n\geq K$.
        Assume that  $\ddagger(10\delta,n)(x,y)$ holds for all 
        pairs of vertices $(x,y)$ in $ B_{v_0} (  R(n)  )$ satisfying   $\star_{10\delta}$.  
 
        Assume furthermore that all parabolic subgroups of $G$ are abelian 
        (so that $X$ satisfies the conclusion of Lemma \ref{lem;horoballs2}).

        For all pairs of points $(x, y)$ in $X$ satisfying $\star_0$, such that 
        $d(x,v_0) \geq R(n)$, the property $\ddagger(0,n)$ is satisfied.

\end{corollary}
\begin{proof}
By  Lemma \ref{lem;move}, if  $x$ and $y$ are in $X_k$, then they satisfy  $\ddagger(0,n)$.  
 On the other hand, if $x$ is not in $X_k$, then it is at depth at least $k$ in an
horoball, and since $y$ is at distance at most $M<k$ from $x$, it is  in the
same horoball.  Lemma  \ref{lem;horoballs2} then 
 guarantees the result.
 \end{proof}

\begin{corollary}\label{coro;carac_ddagger}

   Let $G$ be a relatively hyperbolic group. 
   Assume that the parabolic subgroups are  abelian, and 
    let $X$ be a proper hyperbolic space associated to $G$ as above, and let $K$ be as in paragraph \ref{ss:Constants}.

If  $G$ has connected boundary, and has no proper peripheral splitting, then  there exists $n\geq K$, such that
   $\ddagger(10\delta,n)$ holds for all pairs $(x,y)$  in $ B_{v_0}(  R(n)  )$  satisfying $\star_{10\delta}$.

If  there exists $n\geq K$, such that
   $\ddagger(10\delta,n)$ holds for all pairs $(x,y)$  in $ B_{v_0}(  R(n)  )$  satisfying $\star_{10\delta}$, then $\partial G$ is connected.

\end{corollary}

In fact, what one really needs is  that the parabolic subgroups satisfy the conclusion of Lemma \ref{lem;horoballs2}, in order to apply the previous corollary, and that they are one-ended (or two-ended), 
            finitely presented, and  without infinite torsion subgroup, in order to apply Bowditch's results.

\begin{proof}
If  $G$ has connected boundary, and has no proper peripheral splitting, then by Theorem \ref{theo;bowditch} (1-2), $\partial G$ has no global cut point. Therefore, by  Lemma \ref{lem;BM-cutpoint},  there exists $n$ (which can be assumed to be at least $K$) such that
   $\ddagger(10\delta,n)$ holds for all pairs $(x,y)$  in $ B_{v_0}(  R(n)  )$ satisfying $\star_{10\delta}$.

 Conversely, by Corollary \ref{coro;allpairs}, and Lemma \ref{lem;BM-loc-connected},  $\partial G$ is connected.
 \end{proof}

\section{Algorithms} \label{s:Algorithms}

\subsection{Tools.}

\begin{proposition}\label{prop;check_ddagger}
  There is a procedure that, given a relative presentation of a relatively hyperbolic group $G$ with
  abelian parabolic subgroups (so that the intersection of the generating set with
  each parabolic is sensible), terminates if the boundary $\partial G$ is connected, without global cut point, and only if $\partial G$ is connected. 
\end{proposition}
\begin{proof}
First compute the constant $\delta$ for the space $X$, and set $n=K$. 

Next, compute the ball  $ B_{v_0}(R(n)  +n)$  in $X$. 
  Check whether  $\ddagger(10\delta,n)(x,y)$ holds for all 
        pairs of vertices $(x,y)$ in $ B_{v_0} (  R(n)  )$ satisfying   $\star_{10\delta}$ (there are finitely many such pairs). If this is so, 
then by Corollary \ref{coro;carac_ddagger}, $G$ has connected boundary, so the algorithm is told to stop there. If not, increment $n$ by $1$, and restart this procedure. 

By Corollary \ref{coro;carac_ddagger}, if the procedure does not stop, either the boundary $\partial G$ is disconnected, or it has a global cut point. 
On the other hand, if this procedure {\em does} stop then by Corollary
\ref{coro;carac_ddagger} $\partial X$ is connected.
\end{proof}

\begin{proposition}\label{prop;check_ends}
  There is a procedure that, given a finite presentation of a group, with a solution to the word problem, 
  terminates if and only if the group has more than one end.
\end{proposition}
\begin{proof}
Enumerate the presentations of the group by Tietze transformations. Each time one of them  exhibits a presentation of an amalgamated free product, or of an HNN extension, start in parallel the computation of the multiplication table of the edge group, using the solution to the word problem. One such computation eventually terminates if and only if the edge group is finite. For each of the terminating computations, check whether the splitting is non trivial (if the generators of a vertex group are all in the edge group).  Since the edge group is finite, and the vertex group is given with a finite set of generators, this
involves finitely many applications of the algorithm to solve the word problem. 

The whole procedure is told to terminate when it has found  a non-trivial splitting over a finite group. 
\end{proof}

\begin{proposition}\label{prop;find_splitting}

  There is a procedure that, given a relative presentation of a relatively hyperbolic group $G$, along with a solution to the word problem for $G$, 
  enumerates the peripheral splittings of $G$.

\end{proposition}
\begin{proof}
Enumerate the presentations of the group by Tietze transformations. Each time one of them  exhibits a presentation of 
  a finite bipartite graph of groups,
  start in parallel the following procedure. 

  For each vertex of a given color, of generators  $c_1, \dots c_k$, enumerate the conjugates $c_1^g, \dots, c_k^g$,  for $g\in G$, and check, with a solution to the word problem, whether one of these families coincide with the generators of a given parabolic subgroup. If there is a bijective such correspondence between vertices of one color, and given representative (up to conjugacy) of parabolic subgroups, then the presentation exhibits a peripheral splitting. It is then output in the list.

Given any peripheral splitting, a presentation exhibiting it will appear eventually in the enumeration, and thus will be eventually spotted  by the algorithm. 
\end{proof}

\subsection{Results and Conclusion.}

\begin{theorem} \label{theo;connected_boundary}

There exists an algorithm whose input is a finite presentation of a relatively hyperbolic group with abelian parabolic subgroups, and outputs the answer to the question of the connectivity of its boundary.
\end{theorem}
\begin{proof}
According to Proposition \ref{propB;one-end=connected}, $\partial G$ is disconnected
if and only if $G$ splits nontrivially over a finite group relative to ${\bf P}$.

First of all, we run the procedure from Proposition \ref{prop;check_ends}.  If there
is a nontrivial splitting over a finite subgroup, we will find it.  We can run a procedure
on each such splitting which terminates if each element of ${\bf P}$ is elliptic
({\em i.e.} conjugate into a vertex group).  Thus, if the boundary $\partial G$ is 
disconnected, we will eventually discover this.  We now need a second procedure
to run in parallel which will terminate if the boundary is connected.

Let us run the procedure from Proposition \ref{prop;find_splitting} (recall that by Theorems  \ref{t;dfind}  and  \ref{t:ConetoCusp}, one can compute the relatively hyperbolic structure of the group, a hyperbolicity constant for $X$, and  arbitrarily large balls of $X$, thus one can run \ref{prop;find_splitting}).  The output list of this procedure is a list of graphs of groups (possibly non-terminating).  Start this list with
the trivial graph of groups with one vertex corresponding to $G$ and no edges.

In parallel, every time a new item appears on the list,   let us find  (using the algorithm of \ref{t;dfind}) 
 a relative presentation for the vertex groups, where any adjacent edge groups are parabolic, using the algorithm of \ref{t;dfind}. Then, let us  check whether the  splitting contains a non trivial splitting over a finite group  (the edge groups are given in terms of the presentations of the parabolic subgroups, which are finitely generated abelian, thus it is easy to check whether  they are finite).   
If yes, announce that the group has disconnected boundary (this is true by Proposition \ref{propB;one-end=connected}). In no, then
let us  run the procedure from Proposition \ref{prop;check_ddagger} (checking connectedness of the boundary, modulo presence of cut points) on each of the vertex groups. If the algorithm from Proposition \ref{prop;check_ddagger} has stopped on every non-parabolic vertex group of a given splitting, then announce that the boundary is connected.

We claim that if the procedure from Proposition \ref{prop;check_ddagger} terminates
on each non-parabolic vertex group of a given splitting then the boundary is indeed
connected (as our algorithm announces).  For, in this case  all edge groups are infinite, and  the boundaries of the vertex groups are connected (because the algorithm from Proposition  \ref{prop;check_ddagger} is assumed to have stopped).   One can easily see in this case  that  $\partial G$ is connected, since it is obtained by gluing together boundaries of vertex groups along points fixed by infinite parabolic edge groups, and then compactifying by the boundary of the Bass-Serre tree (see \cite{Dgt} for instance).

It remains to check that  the algorithm presented always terminates. If the boundary is disconnected,
 then by Proposition \ref{propB;one-end=connected} there is a peripheral splitting over a finite group, and it will be listed in the first procedure.  On the other hand if $\partial G$ is connected, by the accessibility Theorem \ref{theoB;access}, the procedure from Proposition \ref{prop;find_splitting} eventually finds a peripheral splitting in which the non parabolic vertex groups have connected boundaries without global cut points. Therefore, by Corollary \ref{coro;carac_ddagger},  on this particular entry at least, the algorithm from Proposition \ref{prop;check_ddagger} terminates. 
\end{proof}

\vskip .5cm

\begin{proof}[Proof of Theorem \ref{t:ExistsSplit}.]
By Theorem  \ref{t;dfind}, one can find a collection of abelian subgroups ${\bf P} = \{P_1, \dots, P_n\}$ (with presentations) of $G$, relative to which it is hyperbolic. One can identify the groups $P_i, i\in I$ that are virtually cyclic, and it is classical that $G$ is still relatively hyperbolic with respect to ${\bf P}'={\bf P} \setminus \{P_i, i\in I\}$, which consists this time of one-ended abelian subgroups.  Now, by Proposition \ref{propB;one-end=connected}, $G$ admits a splitting over a finite group if and only if its boundary (for this latter structure $(G, {\bf P'})$) is connected. Theorem \ref{theo;connected_boundary} can then be applied to $(G, {\bf P'})$, giving the result.
\end{proof}

\begin{proof}[Proof of Theorem \ref{t:Dunwoody}.]
Let us run the procedure of Theorem \ref{t:ExistsSplit}: if the group has one end, there is nothing to add, if it has several ends, the output of the algorithm is a splitting of the group as a graph of groups, with only finite edges groups. The Dunwoody decomposition is a refinement of this splitting, thus we need to run the same procedure on each vertex groups. This process eventually stops, as accessibility is ensured. 
\end{proof}
\begin{proof}[Proof of Theorem \ref{t:Grushko}.]
It suffices to identify the trivial edge groups in the Dunwoody decomposition, that is the output of  Theorem \ref{t:Dunwoody}.  Every edge group in this decomposition is described in the output of \ref{t:Dunwoody} as a subgroup of a parabolic group of $G$ for which we know a presentation, and that is finitely generated abelian. Thus it is easy to decide whether it is trivial. 
\end{proof}

\end{document}